\documentclass[11pt]{article}
\usepackage{geometry}
\usepackage{xcolor}
\usepackage{amsfonts}
\usepackage{mathrsfs}
\usepackage{amsmath}
\usepackage{amsthm}
\usepackage{amssymb}
\def\ds{\displaystyle}
\def\={\buildrel \triangle \over =}
\def\a{\alpha}
\def\b{\beta}

\def\l{\lambda}

\def\t{\times}

\def\th{\theta}
\def\o{\omega}

\def\ns{\noalign{\ss} }

\def\O{\Omega}

\def\ms{\medskip}

\def\q{\quad}
\def\qq{\qquad}

\def\www{\widetilde}

\def\dbR{{\mathbb{R}}}

\def\3n{\negthinspace \negthinspace \negthinspace }
\def\2n{\negthinspace \negthinspace }
\def\1n{\negthinspace }

\def\cA{{\cal A}}

\def\cF{{\cal F}}

\def\pa{\partial}
\def\wt{\widetilde}
\def\cd{\cdot}

\def\div{\hbox{\rm div$\,$}}

\def\|{||}
\def\({\Big (}
\def\){\Big )}
\def\[{\Big[}
\def\]{\Big]}
\def\be{\begin{equation}}
\def\bel{\begin{equation}\label}
\def\ee{\end{equation}}
\def\bt{\begin{theorem}}
\def\bcd{\begin{condition}}
\def\ecd{\end{condition}}
\def\et{\end{theorem}}
\def\bc{\begin{corollary}}
\def\ec{\end{corollary}}
\def\bde{\begin{definition}}
\def\ede{\end{definition}}
\def\bl{\begin{lemma}}
\def\el{\end{lemma}}
\def\bp{\begin{proposition}}
\def\ep{\end{proposition}}
\def\br{\begin{remark}}
\def\er{\end{remark}}
\def\ba{\begin{array}}
\def\ea{\end{array}}
\def\ed{\end{document}}
\def\ns{\noalign{\ms}}
\def\ds{\displaystyle}

\def\om{\omega}

\newtheorem{lemma}{Lemma}[section]
\newtheorem{remark}{Remark}[section]

\newtheorem{theorem}{Theorem}[section]
\newtheorem{corollary}{Corollary}[section]

\newtheorem{definition}{Definition}[section]
\newtheorem{proposition}{Proposition}[section]
\newtheorem{condition}{Condition}[section]

\allowdisplaybreaks

\title{\bf Logarithmic Stability for Coefficients Inverse
Problem of Coupled Wave Equations}
\author{Fangfang Dou\thanks{School of Mathematical Sciences, University of Electronic Science and Technology of China, Chengdu, China. Email: fangfdou@uestc.edu.cn.} \ and Masahiro Yamamoto\thanks{Department of Mathematical Sciences, The University of Tokyo, Komaba, Meguro, Tokyo 153, Japan. Email: myama@ms.u-tokyo.ac.jp.}\ \thanks{Honorary Member of Academy of Romanian Scientists, Splaiul Independentei Street, No. 54, 050094, Bucharest Romania.}\ \thanks{ Peoples’ Friendship University of Russia (RUDN University), 6 Miklukho-Maklaya St, Moscow, 117198, Russian Federation.}}
\date{}                       

\begin{document}

\maketitle
\begin{abstract}
This paper investigates the identification of
two coefficients in a coupled hyperbolic
system with an observation on one component of
the solution. Based on the the Carleman estimate
for coupled wave equations
a logarithmic type stability result is obtained
by measurement data only in a suitably chosen
subdomain under the assumption that the
coefficients are given in a neighborhood of some
subboundary.
\end{abstract}


\section{Introduction and main result}

Let $T>0$ and $\O\subset\dbR^n$ be a nonempty
bounded domain. Write $\mathbf{n}=\mathbf{n}(x)$
for the unit outward normal vector of
$\partial\O$ at $x$. Consider the following
coupled hyperbolic system:
\begin{equation}\label{98eq1}
\left\{
\begin{array}{ll}
\ds \partial_t^2 y_1-\text{div}(a(x)\nabla
y_1)+c_{11}(x)y_1+c_{12}(x)y_2=0 &\text{ in
} Q\triangleq\Omega\times(0,T),\\
\ns\ds\partial_t^2 y_2-\text{div}(a(x)\nabla
y_2)+c_{21}(x)y_1+c_{22}(x)y_2=0 &\text{ in
} Q,\\
\ns\ds \frac{\pa y_1}{\pa\mathbf{n}}=0, \;
\frac{\pa y_2}{\pa\mathbf{n}}=0 & \text{ on
}\pa\O\times(0,T),\\
\ns\ds (y_1(0),\partial_t
y_1(0))=(y_{10},y_{11}), \; (y_2(0),\partial_t
y_2(0))=(y_{20},y_{21})& \text{ in }\Omega.
\end{array}
\right.
\end{equation}

It is well known that wave equations are widely
used to describe many kinds of waves in the
world. In particular, the system \eqref{98eq1}
is a simplified model for describing the
interaction of waves (e.g.,
\cite{Braun1990,Cadoni2007,Kivshar1988}).

Let $\o$ be a nonempty open subset of $\O$. In
this work,  we consider the coefficients inverse
problem for system \eqref{98eq1}, i.e.,
establish the conditional stability for
identifying the coefficients in the zeroth-order
terms $(c_{11},c_{22})$ simultaneously from
suitable observation of one component $y_1$ of
the solution $y=(y_1,y_2)$ in $\o\times(0,T)$.
More precisely, we consider the following
problem:

\ms

\textbf{Problem (IP)}.  Can one recover the
coefficients $(c_{11}, c_{22})$ from a suitable
observation of $y_1$ on $\omega \times(0, T)$?

\ms

Throughout this paper, in order to emphasize the
dependence of the solution of \eqref{98eq1} on
the coefficients, we denote by
$(y_1(c_{11},c_{22}),y_2(c_{11},c_{22}))$ the
solution of \eqref{98eq1} with fixed
coefficients $c_{12}, c_{21}$.

Coefficient inverse problems are important in
various real world applications,  including the
detection and identification of explosives,
nondestructive testing and material
characterization. For their significant
applications, coefficient inverse problems are
widely studied for different equations and
systems. Generally speaking, ``recover" usually
refers to the following two issues:

\begin{itemize}
  \item determining the
coefficients uniquely by the measurement;
  \item giving an
algorithm to compute the coefficients
efficiently.
\end{itemize}
A key step to achieve the above two goals is to
establish an inequality which is called a
stability estimate:
\begin{equation}\label{con-in}
\|(c_{11}-\www{c}_{11},  c_{22}-\www{c}_{22})
\|\leq f(\Vert y_1(c_{11},c_{22}) -
y_1(\www{c}_{11}, \www{c}_{22})\Vert
_{\o\times(0,T)}),
\end{equation}
where $f$ is a non-negative continuous function
satisfying $f(0)=0$.

On one hand, it is clear that if \eqref{con-in}
holds, then $y_1(c_{11},c_{22}) =
y_1(\www{c}_{11}, \www{c}_{22})$ in
$\o\times(0,T)$ implies that $(c_{11}, c_{22}) =
(\www{c}_{11}, \www{c}_{22})$. This implies that
the measurement of $y_1$ in $\o\times(0,T)$ can
uniquely determine the coefficients $(c_{11},
c_{22})$.

On the other hand, according to
\cite{Cheng2000}, one knows that the stability
rate described by the function $f$ is a
quasi-optimal convergence rate of Tikhonov
regularization with a suitable a priori choice
of the regularizing parameters according to
noise levels in data $y_1|_{\o \times(0,T)}$.

In general, there are three common types of $f$:
\begin{enumerate}
  \item $f(\xi)=C\xi$;
  \item $f(\xi)=C\xi^\a$ for some $\a\in (0,1)$;
  \item $f(\xi)=C|\ln\xi|$.
\end{enumerate}
For the first, the second and the third kinds of
$f$,  the estimate \eqref{con-in} indicates
Lipschitz-type stability, H\"older-type
stability and logarithmic-type stability,
respectively.

As one main methodology for the coefficient
inverse problem, we refer to Bukhgeim and
Klibanov \cite{BK}. See also Bellassoued and
Yamamoto \cite{BY2017}, Klibanov \cite{Kli},
Klibanov and Timonov \cite{Klibanov2004} for
example.   The arguments are based on Carleman
estimates, which we discuss. There have been
many works: Beilina, Cristofol, Li and Yamamoto
\cite{BCLY2018}, Bellassoued \cite{B2004},
Bellassoued and Yamamoto \cite{BY2006},
Cannarsa, Floridia and Yamamoto \cite{CFY},
Cannarsa, Floridia, G\"olgeleyen and Yamamoto
\cite{CFGY}, Imanuvilov and Yamamoto
\cite{IY2003}, Klibanov \cite{Klibanov2013},
L\"u and Zhang \cite{Lu2015}, Yu, Liu and
Yamamoto \cite{YLY2018} and the references
therein.  Here we do not intend a comprehensive
list.

Compared with the case of single partial
differential equations, there are much fewer
works addressing coefficients inverse problems
for coupled systems.  By the character of the
Carleman estimate, the inverse problems for
weakly coupling systems, which mean that the
terms of the second order are not coupled, can
be done very similarly to the case of a single
equation if we adopt data of all the components
of the solution. However for strongly coupling
cases, it is more difficult to establish
underlying Carleman estimates and there are very
few researches for inverse problems by Carleman
estimates. As for inverse problems for the
Lam\'e systems which are strongly coupled, see
e.g., Bellassound, Imanuvilov and Yamamoto
\cite{BIY}, Bellassoued and Yamamoto
\cite{BY2007, BY2017}, Imanuvilov and Yamamoto
\cite{IY2005}, for instance.

Our main target is a weakly coupling hyperbolic
system (1), and we describe our main
achievements for the inverse problem:
\begin{itemize}
\item
Data of one component of the solution:\\
As simliar works, we can refer to Benabdallah,
Cristofol, Gaitan and Yamamoto \cite{BCGY},
Alabau-Boussouira, Cannarsa and Yamamoto
\cite{A-BCY} for example.
\item
Data for the inverse problem can be restricted
to an arbitrarily fixed subdomain $\omega$: For
a single wave equation, see \cite{B2004},
\cite{BY2006}. The argument is based on the
Fourier-Bros-Iagolnitzer transform which is a
kind of truncated Laplace transform, and
applications to coupling systems require
non-trivial consideration.
\end{itemize}

In this paper, we establish a logarithmic-type
stability with the measurement on only ONE
component of the solution.  In order to present
the main result, let us  introduce some
notations and conditions. Throughout this paper,
we assume that $a=a(x) \in C^4(\overline{\O})$
satisfying
$$
a>\theta_1 \quad \mbox{on $\overline{\O}$},
\quad \|a\|_{C^4(\overline\O)}\leq M_0, \quad
\Big|\frac{\nabla a(x)\cdot(x-x_0)}{2a(x)}\Big|
\le 1-\theta_0, \ x\in\overline{\O\backslash\o},
$$
from some constants $M_0, \theta_0 > 0$ and $0<
\theta_1\le 1$, and subdomain $\o$ of $\O$.

\begin{remark}
The above assumption on $a$ is for  Lemma
\ref{lm1}. More precisely, it is used to
establish suitable Carleman estimate for
\eqref{98eq1}, which is the key tool to prove
Lemma \ref{lm1}. It is a kind of pseudoconvex
condition and has already been used by several
authors (e.g. \cite{IY2001,Kli}).
\end{remark}
Let $\o_1 \subset \O$ be such that
\begin{equation}\label{tilde o}
\begin{cases}
\ds \o\subset\o_1,\\
\ns\ds\mbox{ for some }x_0\notin \overline\O,\;\{x\in\pa\O: (x-x_0)\cd\nu(x)\geq 0\}\subset\pa\o_1,\\
\ns\ds {\rm
dist}(\pa\o_1\setminus\pa\O,
\pa\o\setminus\pa\O)>0. 
\end{cases}
\end{equation}
Let $\tilde\o$ be a
neighborhood of $\pa\O$ such that
$\o_1\subset\tilde\o$ and ${\rm
dist}(\pa\o_1\setminus\pa\O,
\pa\tilde\o\setminus\pa\O)>0$.

Let us now define the
admissible set of unknown coefficients. Fix
constants  $M_1>0$, $\varpi_1,\varpi_2\in
W^{1,\infty}(\tilde\o)$ and let
$\mathcal{A}=\mathcal{A}(T,\omega,
M_1,\varpi_1,\varpi_2)$ be the set of pairs of
real valued functions $(c_{11},c_{22})$ such
that
\begin{equation}\label{1011eq5}
\begin{array}{ll}\ds
\mathcal{A}=\Big\{(c_{11},c_{22})\in
W^{1,\infty}(\O)^2:\,
\|c_{jj}\|_{W^{1,\infty}(\overline\O)}\leq M_1,
\;    c_{jj} = \varpi_j\text{ for }j=1,2 \Big\}.
\end{array}
\end{equation}

For $s>\frac{3}{2}$, set
$$
X^s(\O)\= \Big\{u\in H^s(\O):\, \frac{\pa u}{\pa
\mathbf{n}}=0\Big\}.
$$

By the classical well-posedness result for wave
equations, similarly to \cite[Lemma
2.1]{IY2001}, for any $(y_{10},y_{11}),
(y_{20},y_{21})\in X^3(\O)\times X^2(\O)$, the
equation \eqref{98eq1} has a unique solution
$$
(y_1,y_2)\in \big[C([0,T];H^3(\O))\times
C^1([0,T];H^2(\O))\times
C^2([0,T];H^1(\O))\big]^2
$$
satisfying that
\begin{equation}\label{12.4-eq1}
\begin{array}{ll}\ds
|(y_1,y_2)|_{[C([0,T];H^3(\O))\cap
C^1([0,T];H^2(\O))\cap
C^2([0,T];H^1(\O))]^2}\\
\ns\ds\leq
C(M_1)(|(y_{10},y_{11})|_{H^3(\O)\cap
H^2(\O)}+|(y_{20},y_{21})|_{H^3(\O)\cap
H^2(\O)}).
\end{array}
\end{equation}
\begin{remark}
The admissible set $\mathcal{A}$ defined by (6),
poses constraints on unknown coefficients:
\begin{itemize}
\item
A priori bounds for $(c_{11},c_{22})$: This is
reasonable because in a physical model, one
usually have some rough estimate on the
coefficients.

\item
We assume to know the values of
$(c_{11}(x),c_{22}(x))$, $x\in\tilde\o$: This
can be interpreted by that one can directly know
physical properties near the boundary.
\end{itemize}
According to the classical well-posedness of
wave equations(e.g., \cite{Ikawa2000}), we know
that there are plenty of solutions such that
$\mathcal{A}$ is nonempty.
\end{remark}

Next, we give the condition for $c_{12}$ and
$c_{21}$:

$\{c_{12},c_{21}\}\subset
W^{2,\infty}(\O)$ and there is a constant
$c_0>0$ such that
\begin{equation}\label{522eq5}
c_{21}\geq c_0 \; \text{ or } -c_{21}\geq c_0
\text{ in } \omega.
\end{equation}
\begin{remark}\label{rem1}
Condition \eqref{522eq5} means that  $y_1$ can
effect $y_2$ adequately. Without \eqref{522eq5},
one cannot obtain information of $y_2$ from
$y_1$ and the observation on $y_1$ is not enough
to determine the coefficients $(c_{11},c_{22})$.
\end{remark}

Now we are ready to state the main result of
this paper.

\begin{theorem}\label{main2D}
There exists $T_0>0$ such that for all $T>T_0$,
we have that
\begin{equation}\label{main2D-eq1}
\begin{array}{ll}\ds
\Vert (c_{11}-\tilde{c}_{11},c_{22}-\tilde{c}_{22})\Vert_{L^2(\Omega)} \\
\ns\ds\leq C\( |\ln \|\partial_t^j(y_1(c_{11},
c_{22}) - y_1(\widetilde{c}_{11},
\widetilde{c}_{22})\|_{L^2(\omega\times (0,T))}|^{-1}C(M_1)\\
\ns\ds\qq + \|\partial_t^j(y_1(c_{11}, c_{22}) -
y_1(\widetilde{c}_{11},
\widetilde{c}_{22})\|_{L^2(\omega\times
(0,T))}\)
\end{array}
\end{equation}
for all $(c_{11},c_{22}), (\tilde c_{11},\tilde
c_{22}) \in \cA$, where $C=C(T) > 0$ is a
constant.
\end{theorem}

From the proof of Theorem \ref{main2D}, one can
see that it can be generalized to a system
coupled by more than two wave equations by data
of reduced numbers of components of data, but in
this paper, we do not pursue the full technical
generality for presenting the key in a simple
way. Following the method in \cite{BCLY2018}, We
can discuss similar inverse problems of
determining all the coefficients but we need to
choose suitable initial values and repeat taking
data. Furthermore we can establish a stability
estimate in determining other combinations such
as $(a, c_{12})$ of two coefficients among $a,
c_{11}, c_{12}, c_{21}, c_{22}$ by a single
measurement of $y_1$ in $\omega \times (0,T)$,
but we do not discuss here. Moreover, the
elliptic operator $\div(a\nabla)$ can be
generalized to a more general one as
$\sum_{j,k=1}^n \pa_{x_j}(a^{jk}\pa_{x_k})$ for
suitable $\{a^{jk}\}_{1\leq j,k\leq n}$. Indeed,
by \cite[Theorems 4.2 and 4.3]{FLZ1}, we can
prove a similar result for Lemma \ref{lm1}. Then
the rest of the proof is similar.
\\

The rest of this paper is organized as follows.
Section 3 is devoted to  presenting   some
auxiliary result, i.e., a H\"older type
stability estimate for the problem we consider
with measurement of one component in an open
subset of the domain satisfying some geometry
conditions, as well as the introduction for the
Fourier-Bros-Iagolnitzer transform. Then in
Section 3, we give the proof of Theorem
\ref{main2D}.




\section{Preliminaries}

In this section, we give some preliminary
results. We first recall a Lipschitz type
estimate.

Let $O_1$ be a subset of $\tilde\o$ such that
$ O_{1}\subset \tilde\o$, $\pa \tilde\om \subset \pa
O_1$.

\begin{lemma}\label{lm1}
For all $T>0$ satisfying
\begin{equation}\label{eq1}
T> \sup_{x\in\O}|x-x_0| \;\mbox{ for the
$x_0$ given in \eqref{tilde o}},
\end{equation}
there exists  
$ O_{1}\subset \tilde\o$ such that $\pa \tilde\om \subset \pa
O_1$ and  a constant $C
>0$ such that for all $(c_{11},c_{22}),$ $(\tilde{c}_{11},\tilde{c}_{22}) \in \cA$,
\begin{equation}\label{1009eq1}
\begin{array}{ll}\ds
\|(
c_{11}-\tilde{c}_{11},c_{22}-\tilde{c}_{22})\|_{L^2(\O)}\\
\ns\ds\leq C
\sum_{j=1}^2\sum_{k=1}^2\|\pa_t^{k}(y_j(c_{11},c_{22})-
\tilde{y}_j(\tilde{c}_{11},\tilde{c}_{22}))
\|_{L^2((0,T)\t(\o_1\setminus O_1))}.
\end{array}
\end{equation}
\end{lemma}

Lemma \ref{lm1} can be obtained directly by
following the proof of Theorem 1.1 in
\cite{IY2001} step by step. We omit it here.

Next, we give a brief introduction for
Fourier-Bros-Iagolnitzer transform (shortened to
F.B.I. transform), which is a generalization of
the Fourier transform in this subsection. A
detailed introduction to this can be found in
\cite{Delort1992}. As an important application
to a hyperbolic equation, see Robbiano
\cite{Ro}, and here we modify the aruments in
\cite{Ro}. Let
$$
F(z)\triangleq\frac1{2\pi}\int_\mathbb{R}
e^{iz\tau}e^{-\tau^2}d\tau.
$$
Then
\begin{equation*}
F(z)=\frac{\sqrt{\pi}}{2\pi}e^{\frac14
(|\text{Im} z|^2-|\text{Re} z|^2)}e^{-\frac i2
(\text{Im} z \text{Re} z)}.
\end{equation*}
For every $\lambda\geq1$, define
\begin{equation*}
F_\lambda(z)\triangleq\lambda F(\lambda
z)=\frac{1}{2\pi}\int_\mathbb{R}
e^{iz\tau}e^{-(\frac{\tau}{\lambda})^2} d\tau.
\end{equation*}
It can be easily seen that
\begin{equation*}
|F_\lambda(z)|=\frac{\sqrt{\pi}}{2\pi}\lambda
e^{\frac{\lambda^2}4 (|\text{Im} z|^2-|\text{Re}
z|^2)}.
\end{equation*}

Let $s, l_0\in\mathbb{R}$ and recall that
$i=\sqrt{-1}$. The F.B.I. transformation
$\mathcal{F}_\lambda$ for
$f\in\mathcal{S}(\mathbb{R}^{n+1})$ is defined
as follows:
\begin{equation}\label{eq19}
(\cF_\l f)(x,s)\triangleq\int_{\mathbb{R}}
F_\lambda(l_0+is-l)\Phi(l)f(x,l)dl.
\end{equation}
%


\section{Proof of Theorem \ref{main2D}}


Let $\omega_0$ is an arbitrary fixed
nonempty subset of $ \omega$ such that
$\overline{\omega_0}\subset \omega$.  By \cite[Lemma 1.1]{FI1996}, we
know that there exists a function $\hat\psi\in
C^2(\tilde{\omega})$ such that
\begin{equation}\label{53eq1}
\begin{cases}\ds
\hat{\psi}(x)>0, \q  \forall x\in \tilde\omega, \\
\ns\ds \hat{\psi}(x)=0,\q \forall x\in
\partial \tilde\omega, \\
\ns\ds |\nabla\hat{\psi}(x)|>0, \q  \forall
x\in\overline{\tilde\omega\backslash \o_0}.
\end{cases}
\end{equation}
We can conclude from \eqref{53eq1} that there
exist two constants $\beta_1,\beta_2>0$, where
$\beta_2>(\beta_1+\|\hat\psi\|_{L^\infty(\tilde\o)})/2$,
and $(\o_1\setminus O_1)\subset \subset O_2\subset \subset O_3\subset \subset \tilde\o$ such that
\begin{equation}\label{520eq1}
\hat{\psi}(x)\leq\beta_1, \q \forall x\in O_3\setminus O_2
\end{equation}
and that
\begin{equation}\label{520eq3}
\hat{\psi}(x)\geq\beta_2, \q \forall x\in
\omega_1\setminus O_1.
\end{equation}
It follows from the last condition in
\eqref{53eq1} that the maximum value of
$\hat\psi$ can only be attained in $\omega_0$,
i.e., there exists a point $\hat x\in\omega_0$
such that
\begin{equation}\label{520eq2}
\hat{\psi}(\hat x)=\max_{x\in
\omega}\hat\psi(x).
\end{equation}

Let $\chi\in C^\infty(\tilde\o)$  satisfying that
$0\leq\chi\leq1$ and
\begin{equation}\label{eq20}
\chi(x)=\left\{
\begin{array}{ll}
1, &\text{ if } x\in  O_2,\\
0, &\text{ if } x\in \tilde\o\setminus O_3.
\end{array}
\right.
\end{equation}
Put
$$
\begin{array}{ll}\ds
w_j=(y_j(c_{11},c_{22})-\tilde{y}_j(\tilde{c}_{11},\tilde{c}_{22})),
\q g_j=(\tilde{c}_{jj}-c_{jj})\chi \tilde y_j.
\end{array}
$$

\ms

Let $(U_1,U_2)=(\chi\partial_t w_1,\chi
\partial_t w_2)$ and $(V_1,V_2)=(\chi
\partial_t^2w_1,\chi \partial_t^2w_2)$,
respectively. By the assumption of $\mathcal{A}$
in \eqref{1011eq5}, there exists $C(M_0,M_1)>0$
such that
\begin{equation}\label{57eq10}
\|(U_1,U_2)\|^2_{L^2(-T,T;H^2(\tilde\omega))^2}
+\|(\pa_tU_1,\pa_tU_2)\|^2_{L^2(-T,T;H_0^1(\tilde\omega))^2}\leq
C(M_0,M_1)^2.
\end{equation}
Then by fundamental calculation, we have that
\begin{eqnarray}\label{eq21}
\left\{\!\!\!
\begin{array}{ll}\ds
\partial_t^2 U_1\!-\!\text{div}(a(x)\nabla U_1)\!  +\!c_{11}U_1\! +\! c_{12}U_2\!=\!-(\nabla a\cdot\nabla \chi) U_1\!-\!a[\Delta, \chi]U_1 &\mbox{in }  \tilde\omega\times(-T,T),\\
\ns\ds \partial_t^2 U_2\!-\!\text{div}(a(x)\nabla U_2) \! +\!c_{21}U_1 \!+ \!c_{22}U_2\!=\!-(\nabla a\cdot\nabla \chi) U_2\!-\!a[\Delta, \chi]U_2 &\mbox{in }  \tilde\omega\times(-T,T),\\
\ns\ds U_1=U_2=0 & \text{in }\partial \tilde\omega\t(-T,T),\\
\ns\ds (U_1(0),\partial_t U_1(0))=(0,\partial_{t}^2 g_1 (x,0))& \text{in } \tilde\o, \\
\ns\ds (U_2(0),\partial_t
U_2(0))=(0,\partial_{t}^2 g_2 (x,0))& \text{in }
\tilde\o.
\end{array}
\right.
\end{eqnarray}
and
\begin{eqnarray}\label{18eq4}
\left\{\!\!\!
\begin{array}{ll}\ds
\partial_t^2 V_1\!-\!\text{div}(a(x)\nabla V_1)\! +\!c_{11}V_1 \!+\! c_{12}V_2\! =\!-(\nabla a\cdot\nabla \chi) v_1\!-\!a[\Delta, \chi]v_1 & \text{in } \tilde\o\times(-T,T),\\
\ns\ds \partial_t^2 V_2\! -\!\text{div}(a(x)\nabla V_2)\!+\!c_{21}V_1 \!+\! c_{22}V_2 \!=\!-(\nabla a\cdot\nabla \chi) v_2\!-\!a[\Delta, \chi]v_2  & \text{in } \tilde\o\times(-T,T),\\
\ns\ds V_1=0, V_2=0 &  \text{on
}\partial \tilde\o\times(-T,T),\\
\ns\ds (V_1(0),\partial_t V_1(0))=(\partial_{t}^2 g_1 (x,0),\partial_t^3 g_1 (x,0))& \text{in } \tilde\o, \\
\ns\ds (V_2(0),\partial_t V_2(0))=(\partial_{t}^2
g_2 (x,0),\partial_t^3g_2 (x,0))& \text{in } \tilde\o.
\end{array}
\right.
\end{eqnarray}

Let $\Phi\in C^\infty_0(\mathbb{R})$ satisfying
the following conditions:
$$
\begin{cases}\ds
\Phi\in C^\infty_0\(\[-\frac L2,\frac
L2\];[0,1]\),\\
\ns\ds \Phi=1 \mbox{ on }\[-\frac L4,
\frac{L}4\],\\
\ns\ds
|\Phi'|\leq\frac1{4L},\q|\Phi''|\leq\frac1{4L},
\end{cases}
$$
where $L>0$ will be chosen later.

Take
$$
K=\[-\frac L2,-\frac
L4\]\bigcup\[\frac{L}4,\frac L2\],\qq
K_0=\[-\frac{L}8,\frac{L}8\].
$$
and let $l_0 \in K_0$ in \eqref{eq19}.

Let $1<b_0<b\leq 2$ and $\xi\in C^\infty(\dbR)$ satisfy
$0\leq\xi\leq 1$ and
\begin{equation}\label{1011eq2}
\xi(t)=\left\{
\begin{array}{ll}
\ds 1, \qq \text{in }[-b_0,b_0],\\
\ns\ds 0, \qq \text{in }(-\infty,-b]\cup[b,+\infty).
\end{array}
\right.
\end{equation}

Let $\wt U_j^F(x,s)$ and $\wt V_j^F(x,s)$ be the F.B.I.
transform of $U_j$ and $V_j$, respectively,
$j=1,2$, i.e.,
$$
\ds\wt U_j^F(x,s)= \int_{\mathbb{R}}
F_\lambda(l_0+is-l)\Phi(l)U_j(x,l)dl
$$
and
$$
\ds \wt V_j^F(x,s)= \int_{\mathbb{R}}
F_\lambda(l_0+is-l)\Phi(l)V_j(x,l)dl.
$$
Let $U_j^F(x,s)=\xi(s)\wt U_j^F(x,s)$ and $V_j^F(x,s)=\xi(s)\wt V_j^F(x,s)$.

Since
$$
\begin{array}{ll}\ds
\q \partial_s^2 \wt U^F_j(x,s) = -\xi(s)\partial_s\int_{\mathbb{R}} i\partial_l F_\lambda(l_0+is-l)\Phi(l)U_j(x,l)dl\nonumber\\
\ns\ds=  i \partial_s \int_{\mathbb{R}}
F_\lambda(l_0+is-l)\left(\Phi'(l)U_j(x,l)+\Phi(l)\partial_l U_j(x,l)\right)dl\nonumber\\
\ns\ds= -\int_{\mathbb{R}}
F_\lambda(l_0+is-l)\left(\Phi''(l)U_j(x,l)+2\Phi'(l)\partial_l
U_j(x,l)+\Phi(l)\partial_l^2 U_j(x,l)\right)dl,
\end{array}
$$
we have that
\begin{equation}\label{911eq1}\left\{
\begin{array}{ll}\ds
\partial_s^2 U^F_1 +\text{div}(a(x)\nabla U^F_1) -c_{11}U^F_1-c_{12}U^F_2 =G_1+H_1+L_1  &\mbox{ in } \tilde\omega\times\mathbb{R},\\
\ns\ds \partial_s U^F_2 +\text{div}(a(x)\nabla U^F_2)  -c_{21}U^F_1-c_{22}U^F_2 =G_2 +H_2+L_2 &\mbox{ in } \tilde\omega\times\mathbb{R},\\
\ns\ds \frac{\pa U^F_1}{\pa\nu}=\frac{\pa
U^F_2}{\pa\nu}=0 &\mbox{ on }\partial
\tilde{\omega}\times\mathbb{R},
\end{array}\right.
\end{equation}
where
$$
\begin{array}{ll}\ds
 G_j(x,s) =\int_{\mathbb{R}}
F_\lambda(l_0+is-l)\Phi(l)(\nabla a\cdot\nabla
\chi) U_j+a[\Delta, \chi]U_j)dl,\\
\ns\ds H_j(x,s) =-\int_{\mathbb{R}}
F_\lambda(l_0+is-l)(\Phi''(l)U_j(x,l)+2\Phi'(l)\partial_l
U_j(x,l))dl,\\
\ns\ds L_j= 2\xi_s\pa_s\wt U^F_j(x,s) + \xi_{ss}\wt U^F_j(x,s).
\end{array}
$$

Let

\begin{equation}\label{12.4-eq5}
\th=e^{\ell}, \q \ell=\zeta\phi, \q
\phi=e^{\mu\psi}, \q
\psi=\psi(x,s)\triangleq\frac{\hat\psi(x)}{M\|\hat\psi\|_{L^\infty(\tilde\o)}}+b^2-s^2.
\end{equation}
%
Here  
$$\frac{1-\frac{\b_2}{\|\hat\psi\|_{L^\infty(\tilde\o)}}}{b_0^2-1}<M<\frac{1-\frac{\b_1}{\|\hat\psi\|_{L^\infty(\tilde\o)}}}{b_0^2},$$
where $\lambda$ and $\mu$ are parameters,  and $s\in (-b,b)$.

By (4.33) in \cite{Lu2013}, there exists a
constant $\mu_0>0$ such that for all
$\mu\geq\mu_0$, one can find two constants
$C=C(\mu)>0$ and $\zeta_0=\zeta_0(\mu)$ so that
for all $\zeta\geq \zeta_0$,  the solution
$(U^F_1,U^F_2)\in H^1((-b,b)\t\tilde\o)^2$ to
\eqref{911eq1} satisfies that
\begin{eqnarray}\label{12.4-eq2}
&& \3n\3n\zeta\mu^2\int_{-b_0}^{b_0} \int_{
{\tilde\omega}}
\th^2\phi\big(|\nabla U^F_1|^2\!+\!|\partial_s U^F_1|^2\!+\!\!\zeta^2\mu^2\phi^2| U^F_1|^2\!+\!|\nabla U^F_2|^2\!+\!|\partial_s U^F_2|^2\!+\!\zeta^2\mu^2\phi^2| U^F_2|^2\big)dxds\nonumber\\
&& \3n\3n
\leq C \[\int_{-b}^b\int_{ {\tilde\omega}}\th^2\big(|G_1(x,s)+H_1(x,s)|^2+|G_2(x,s)+H_2(x,s)|^2\big)dxds\nonumber\\
&& \3n\3n\q\q +\zeta \mu^2\int_{-b}^b\int_{ \o_0}\th^2\phi (|\nabla U^F_1|^2\!+\!|\partial_s U^F_1|^2+|\nabla U^F_2|^2\!+\!|\partial_s U^F_2|^2)dxds\\
&& \3n\3n\q\q  +\zeta^3\mu^4\int_{-b}^b\int_{ \o_0}\th^2\phi^3(| U^F_1|^2 + | U^F_2|^2)dxds\]\nonumber\\
&& \3n\3n\q
+C\int_{(-b,-b_0)\cup(b_0,b)}\int_{\tilde\o}\th^2\big(|\partial_s
\wt U^F_1|^2+|\wt U^F_1|^2+|\partial_s\wt U^F_2|^2+|
\wt U^F_2|^2\big)dxds.\nonumber
\end{eqnarray}
Let us get rid of the terms of $U^F_2$ in the
second and third integrals in the right hand
side of \eqref{12.4-eq2}.

Let $\o_{0,j}$ ($j=1,2$) satisfy that
$\o_0\subset\subset
\o_{0,1}\subset\subset\o_{0,2}\subset\subset
 \o$. We choose cutoff
functions $\eta_{j} \in
C^{\infty}(\omega;[0,1])$ $(j=1,2,3)$
satisfying
\begin{equation}\label{12.4-eq3}
\left\{\begin{array}{ll} \eta_{j}(x)=1, &
\forall x \in \o_{0,j-1}, \\
\ns\ds 0<\eta_{j} \leq 1, & \forall x \in
\omega_{0,j}, \\
\ns\ds  \eta_{j}(x)=0, &  \forall x \in \omega
\backslash \omega_{0,j}.
\end{array}
\right.
\end{equation}
It is easy to see that
\begin{equation}\label{12.4-eq4}
\begin{array}{ll}\ds
\theta^{2} \phi  \eta_{1}^{2}
\overline{U^F_2}\big[\partial_{s}^2
U^F_{2}+\text{div}(a(x)\nabla U^F_2)\big]
\\ \ns\ds =\partial_s \left(\theta^{2} \phi \eta_{1}^{2}
\overline{U^F_2} \partial_s U^F_{2}\right)
-\theta^{2} \phi \eta_{1}^{2}\left|\partial_s
U^F_{2}\right|^{2}-\left(\theta^{2} \phi
\eta_{1}^{2}\right)_{s} \overline{U^F_2}
\partial_s U^F_{2}+\div\left(\theta^{2} \phi
\eta_{1}^{2} a \overline{U^F_2} \nabla U^F_2\right) \\
\ns\ds \q -\theta^{2} \phi \eta_{1}^{2} a|\nabla
U^F_2|^2- a \nabla\left(\theta^{2} \phi
\eta_{1}^{2}\right) \overline{U^F_2}\nabla
U^F_2.
\end{array}
\end{equation}
Integrating \eqref{12.4-eq4} on $(-b, b) \times
\tilde\omega$ and noting that $U^F_2(-b)=U^F_2(b)=0$
in $\tilde\omega$, by \eqref{911eq1} and
\eqref{12.4-eq5},  we see that there exists
$\zeta_1>0$ such that for all $\zeta \geq
\zeta_{1}$,
\begin{equation}\label{12.4-eq8}
\begin{array}{ll}
\ds  \int_{-b}^{b} \int_{\omega_{0}} \theta^{2}
\phi\big( |\nabla U^F_2|^{2}+\left|\partial_s
U^F_{2}\right|^{2}\big)
d x d s  \\
\ns\ds   \leq C\[\int_{-b}^b\int_{
\tilde{\omega}}\th^2 |G_2(x,s)+H_2(x,s)|^2
dxds+\zeta^{2} \mu^{2} \int_{-b}^{b}
\int_{\omega_{0,1}}^{b} \theta^{2}
\phi^{3}|U^F_2|^{2} d x d s \\
\ns\ds\q + \int_{(-b,-b_0)\cup(b_0,b)}\int_{\tilde\o}\th^2\big(|\partial_s
\wt U^F_2|^2+|\wt U^F_2|^2 \big)dxds \].
\end{array}
\end{equation}
Now we estimate $\ds\int_{-b}^{b}
\int_{\omega_{0,1}} \theta^{2}
\phi^{3}|U^F_2|^{2} d x d s $. It is easy to see
that
\begin{equation}\label{12.4-eq6}
\begin{array}{ll}\ds
\theta^{2} \phi^{3} \eta_{2}^{3}
\overline{U^F_2}\left[\partial_{s}^2
U^F_{1}+\div(a
\nabla U^F_1) \right]\\
\ns\ds =\theta^{2} \phi^{3} \eta_{2}^{3}
U^F_1\left[\partial_{s}^2
\overline{U^F_{2}}+\div\left(a \nabla
\overline{U^F_2}\right) \right]+\partial_s
\left[\theta^{2} \phi^{3}
\eta_{2}^{3}\left(\overline{U^F_2}
\partial_s U^F_{1}-\partial_s \overline{U^F_{2}}
U^F_{1}\right)\right] \\
\ns\ds\q -\partial_s \left(\theta^{2} \phi^{3}
\eta_{2}^{3}\right) \overline{U^F_2}
\partial_s U^F_{1}+\partial_s \left[\partial_s \left(\theta^{2} \phi^{3}
\eta_{2}^{3}\right) \overline{U^F_2}
U^F_1\right] -\partial_s \left[\partial_s
\left(\theta^{2} \phi^{3} \eta_{2}^{3}\right)
U^F_1\right] \overline{U^F_2}\\
\ns\ds\q +\div\left[\theta^{2} \phi^{3}
\eta_{2}^{3} a \left(\overline{U^F_2} \nabla
U^F_1-U^F_1\nabla\overline{U^F_2} \right)\right]
-  a \nabla\left(\theta^{2} \phi^{3}
\eta_{2}^{3}\right)\overline{U^F_2}
\nabla U^F_1  \\
\ns\ds\q +\div\left[a \nabla\left(\theta^{2}
\phi^{3} \eta_{2}^{3}\right)\overline{U^F_2}
U^F_1\right] -\div\left[a \nabla\left(\theta^{2}
\phi^{3} \eta_{2}^{3}\right) U^F_1\right]
\overline{U^F_2}.
\end{array}
\end{equation}
Integrating \eqref{12.4-eq6} on $(-b,
b)\times\tilde\o$ and noting that $U^F_1(-b) =
U^F_1(b) = 0$ in $\tilde\o$, by \eqref{911eq1}, we
find that
\begin{equation}\label{12.4-eq7}
\begin{array}{ll}\ds
\int_{-b}^{b} \int_{\omega_{0,1}}  \theta^{2}
\phi^{3}|U^F_2|^{2} d x d s\\
\ns\ds\leq C\[\int_{-b}^{b} \int_{\tilde\omega}
\th^2\big(|G_1(x,s)+H_1(x,s)|^2+|G_2(x,s)+H_2(x,s)|^2\big)dxds
\\
\ns\ds\qq + \zeta^2\mu^2\int_{-b}^{b}
\int_{\omega_{0,2}}\th^2\phi^5(|\nabla
U^F_1|^{2} + |\nabla
U^F_{1,s}|^{2}+\mu|U^F_1|^{2})dxds\\
\ns\ds\qq + \int_{(-b,-b_0)\cup(b_0,b)}\int_{\tilde\o}\th^2\big(|\partial_s
\wt U^F_1|^2+|\wt U^F_1|^2+|\partial_s\wt U^F_2|^2+|
\wt U^F_2|^2\big)dxds 
\].
\end{array}
\end{equation}

Similar to \eqref{12.4-eq8}, we can obtain that
\begin{equation}\label{12.4-eq9}
\begin{array}{ll}\ds
\int_{-b}^{b}
\int_{\omega_{0,2}}\th^2\phi^5(|\nabla
U^F_1|^{2} + |\nabla U^F_{1,s}|^{2} )dxds \\
\ns\ds   \leq C\[\int_{-b}^b\int_{
\tilde{\omega}}\th^2 |G_1(x,s)+H_1(x,s)|^2
dxds+\zeta^{2} \mu^{2} \int_{-b}^{b}
\int_{\omega_{0}}^{b} \theta^{2}
\phi^{7}|U^F_1|^{2} d x d s\\
\ns\ds\qq + \int_{(-b,-b_0)\cup(b_0,b)}\int_{\tilde\o}\th^2\big(|\partial_s
\wt U^F_1|^2+|\wt U^F_1|^2\big)dxds\].
\end{array}
\end{equation}

Combing \eqref{12.4-eq2}, \eqref{12.4-eq8},
\eqref{12.4-eq7} and  \eqref{12.4-eq9}, we know
 there exists a constant $\mu_2>0$ such that for
all $\mu\geq\mu_2$, one can find two constants
$C=C(\mu)>0$ and $\zeta_2=\zeta_2(\mu)$ so that
for all $\zeta\geq \zeta_2$,  the solution
$(U^F_1,U^F_2)\in H^1((-b,b)\t\o)^2$ to
\eqref{911eq1} satisfies that
\begin{eqnarray}\label{911eq2}
&& \3n\3n\l\mu^2\int_{-b_0}^{b_0} \int_{
\tilde{\omega}}
\th^2\phi\big(|\nabla U^F_1|^2\!+\!|\partial_s U^F_1|^2\!+\!\!\l^2\mu^2\phi^2| U^F_1|^2\!+\!|\nabla U^F_2|^2\!+\!|\partial_s U^F_2|^2\!+\!\l^2\mu^2\phi^2| U^F_2|^2\big)dxds\nonumber\\
&& \3n\3n
\leq C \[\int_{-b}^b\int_{\tilde{\omega}}\th^2\big(|G_1(x,s)+H_1(x,s)|^2+|G_2(x,s)+H_2(x,s)|^2\big)dxds\nonumber\\
&& \3n\3n\q\q\q \q\q +\int_{-b}^b\int_{\o_0}\l^5\mu^7\th^2\phi^5| U^F_1|^2dxds\]\\
&& \3n\3n\q 
+C\int_{(-b,-b_0)\cup(b_0,b)}\int_{\tilde\o}\th^2\big(|\partial_s
\wt U^F_1|^2+|\wt U^F_1|^2+|\partial_s \wt U^F_2|^2+|
\wt U^F_2|^2\big)dxds.\nonumber
\end{eqnarray}

\ms

Set $ \phi_j=e^{\mu\psi_j}, j=1,2,3,4 $, where

\begin{eqnarray}
&&\psi_1=\frac{\beta_1}{M\|\hat{\psi}\|_{L^\infty({\tilde\o})}}+b^2, \q \psi_2=\frac{\hat\psi(\hat x)}{M\|\hat{\psi}\|_{L^\infty({\tilde\o})}}+b^2=\frac1M+b^2, \nonumber\\
 &&\psi_3=\frac{\beta_2}{M\|\hat{\psi}\|_{L^\infty({\tilde\o})}}+b^2-1, \q \psi_4=\frac{\hat\psi(\hat x)}{M\|\hat{\psi}\|_{L^\infty({\tilde\o})}}+b^2-b_0^2=\frac1M+b^2-b_0^2.
\end{eqnarray}
%
From the bound for $M$ we know $\psi_1<\psi_4$.
By the property of F.B.I. transformation, we
have that
\begin{eqnarray}\label{eq23}
&&\q\int_{-b}^b\int_{\o_0}\zeta^5\mu^7\th^2\phi^5| U^F_1|^2dxds\nonumber\\
&&\leq \zeta^5\mu^7\max_{x\in\omega_0,s\in[{-b},b]} \big(\th^2\phi^5\big) \int_{-b}^b\int_{\omega_0}|\wt U^F_1(x,s)|^2 dxds\nonumber\\
&&\leq  \zeta^5\mu^7\phi_2^5e^{2\zeta\phi_2} \int_{-b}^b\int_{\omega_0}\Big|\int_{\mathbb{R}} F_\lambda(l_0+is-l)\Phi(l)U_1(x,l)dl \Big|^2 dxds\\
&&\leq   \zeta^5\mu^7\phi_2^5e^{2\zeta\phi_2} \int_{-b}^b\int_{\omega_0}\Big|\int_{\mathbb{R}} \frac{\sqrt\pi}{2\pi} \lambda e^{\frac{\lambda^2}4(s^2-|l_0-l|^2)}\Phi(l)U_1(x,l)dl \Big|^2  dxds\nonumber\\
&&\leq \frac{\lambda^2}{4\pi} \zeta^5\mu^7\phi_2^5e^{2\zeta\phi_2}  \int_{-b}^b e^{\frac{\lambda^2}2 s^2}ds \big|\sup \Phi\big|^2 \int_{\omega_0}\Big|\int_{-\frac L2}^{\frac L2} U_1(x,l)dl \Big|^2  dx\nonumber\\
&&\leq \frac{\lambda^2Lb}{2\pi}
\zeta^5\mu^7\phi_2^5e^{2\zeta\phi_2}
e^{\frac{\lambda^2}2 b^2}
\int_{\omega_0}\int_{-\frac L2}^{\frac L2}
|U_1(x,l)|^2dldx.\nonumber
\end{eqnarray}
From the definition of $H_j$, we see that
\begin{eqnarray}\label{eq17}
&&\q\int_{-b}^b\int_{\tilde\omega} |H_j(x,s)|^2 dxds\nonumber\\
&&\leq \int_{-b}^b\int_{\tilde\omega} \Big|-\int_{\mathbb{R}} F_\lambda(l_0+is-l)(\Phi''(l)U(x,l)+2\Phi'(l)\partial_l U(x,l))dl\Big|^2 dxds\nonumber\\
&&\leq\int_{-b}^b\int_{\tilde\omega} \Big|\int_{K} \frac{\sqrt{\pi}}{2\pi}\lambda e^{\frac{\lambda^2}4 (s^2-|l_0-l|^2)}(\Phi''(l)U(x,l)+2\Phi'(l)\partial_l U(x,l))dl\Big|^2 dxds \\
&&\leq \frac1{2\pi}\lambda^2 e^{\frac{\lambda^2}2 b^2}b\int_{ \tilde\omega} \Big| \int_{K} e^{-\frac{\lambda^2}4 |l_0-l|^2} (\Phi''(l)U(x,l)+2\Phi'(l)\partial_l U(x,l))dl\Big|^2 dx\nonumber\\
&&\leq \frac{b}{2\pi}\lambda^2 e^{\frac{\lambda^2}2 b^2}\frac L2\int_{\tilde\omega} \int_{K} e^{-\frac{\lambda^2}2 |l_0-l|^2} \(2\Big|\Phi''(l)U(x,l)\Big|^2+8\Big|\Phi'(l)\partial_l U(x,l))\Big|^2\)dl dx\nonumber\\
&&\leq \frac{2\lambda^2bL}{\pi}
e^{\frac{\lambda^2}2 (b^2-(\frac L8)^2)}
\max_{K}\{ |\Phi''(l)|^2, |\Phi'(l)|^2\} \int_{
\tilde\omega} \int_{K} \(|U_j(x,l)|^2+|\partial_l
U_j(x,l)|^2\) dldx.\nonumber
\end{eqnarray}

Since $\chi=1$ in $\omega_1$ and
$\chi=0$ in $\tilde\o\setminus O_1$, $G_j(x,s)=0$ in $\o_1\cup(\tilde\o\setminus O_1)$ and therefore
\begin{eqnarray}\label{eq16}
&&\3n\3n\3n\q\int_{-b}^b\int_{\tilde\omega}
|G_j(x,s)|^2
dxds\nonumber\\
&&\3n\3n\3n=\int_{-b}^b\int_{\tilde\omega} \Big|\int_{\mathbb{R}} F_\lambda(l_0+is-l)\Phi(l)((\nabla a \cdot\nabla \chi) U_j+a [\Delta, \chi]U_j)dl\Big|^2 dxds\nonumber\\
&&\3n\3n\3n\leq\int_{-b}^b\int_{\tilde \omega} \Big| \int_{-\frac L2}^{\frac L2}\frac{\sqrt{\pi}}{2\pi}\lambda e^{\frac{\lambda^2}4 (s^2-|l_0-l|^2)} ((\nabla a \cdot\nabla \chi) U_j+a [\Delta, \chi]U_j)dl\Big|^2 dxds\\
&&\3n\3n\3n \leq\frac1{2\pi}\lambda^2
e^{\frac{\lambda^2}2 b^2}b \int_{O_3\setminus
O_2} \Big| \int_{-\frac L2}^{\frac L2}  e^{-\frac{\lambda^2}4 |l_0-l|^2} ((\nabla a \cdot\nabla \chi) U_j+a [\Delta, \chi]U_j)dl\Big|^2 dx\nonumber\\
&& \3n\3n\3n\leq\frac{\lambda^2b L}{2\pi}
e^{\frac{\lambda^2}2
b^2}\|a\|_{C^1(\O)}^2\max\{|\nabla\chi|^2,
|\Delta\chi|^2\} \int_{-\frac L2}^{\frac L2}
\int_{O_3\setminus
O_2} \left( |U_j(x,l)|^2 +
|\nabla
U_j(x,l)|^2\right)dxdl\nonumber\\
&&\3n\3n\3n\leq \frac{\lambda^2b LM_0^2}{2\pi}
e^{\frac{\lambda^2}2 b^2}\max\{|\nabla\chi|^2,
|\Delta\chi|^2\} \int_{-\frac L2}^{\frac L2}
\int_{O_3\setminus
O_2} \left( |U_j(x,l)|^2 +
|\nabla U_j(x,l)|^2\right)dxdl.\nonumber
\end{eqnarray}

Consequently,
\begin{equation}\label{921eq1}
\!\!\begin{array}{ll}
\ds\q\int_{-b}^b\int_{\wt\omega} \th^2|H_j(x,s)|^2  dxds\\
\ns\ds\leq\!
e^{2\zeta\phi_2}\frac{2\lambda^2 bL}{\pi}
e^{\frac{\lambda^2}2 (b^2-(\!\frac L8\!)^2)}
\max_{K}\{ |\Phi''(l)|^2, |\Phi'(l)|^2\}\!
\int_{\wt\omega}\! \int_{K}\!
\(|U_j(x,l)|^2\!+\!|\partial_t U_j(x,l)|^2\!\)
dldx,
\end{array}
\end{equation}
\begin{equation}\label{92eq2}
\begin{array}{ll}
\ds\q\int_{-b}^b\int_{\wt\omega} \th^2|G_j(x,s)|^2  dxds\\
\ns\ds\leq e^{2\zeta\phi_1}\frac{\lambda^2b
LM_0^2}{2\pi} e^{\frac{\lambda^2}2
b^2}\max\{|\nabla\chi|^2, |\Delta\chi|^2\}
\int_{-\frac L2}^{\frac L2} \int_{O_3\setminus
O_2} \left( |U_j(x,l)|^2 + |\nabla
U_j(x,l)|^2\right)dxdl.
\end{array}
\end{equation}

For simplicity of notations, without loss of
generality, we assume that $T=1$. Substituting
\eqref{eq23},  \eqref{921eq1}, and \eqref{92eq2}
into \eqref{911eq2}, we obtain that
\begin{eqnarray}\label{911eq3}
&&\3n\3n \3n\ds
\zeta^3\mu^4\phi_3^3e^{2\zeta\phi_3}\int_{-1}^{1}
\int_{\tilde\omega }
\(| U^F_1|^2+| U^F_2|^2\)dxds\nonumber\\
&&\3n\3n\3n\leq \zeta^3\mu^4\int_{-b_0}^{b_0}
\int_{\tilde\omega}
\th^2\phi^3\(| U^F_1|^2+| U^F_2|^2\)dxds\nonumber\\
&&\3n\3n\3n \leq\zeta\mu^2\int_{-b_0}^{b_0}\!
\int_{\tilde{\omega}} \th^2\phi\big(|\nabla
U^F_1|^2\!+\!|\partial_s U^F_1|^2\!+
\!\zeta^2\mu^2\phi^2| U^F_1|^2\!+\!|\nabla U^F_2|^2\!+\!|\partial_s U^F_2|^2\!+\!\zeta^2\mu^2\phi^2| U^F_2|^2\big)dxds\nonumber\\
&&\3n\3n\3n
\leq C \[e^{2\zeta\phi_2}\frac{4\lambda^2bL}{\pi} e^{\frac{\lambda^2}2 (b^2-(\frac L8)^2)} \max_{K}\{ |\Phi''(l)|^2, |\Phi'(l)|^2\} \int_{\tilde\omega} \int_{K} \(|U_1(x,l)|^2+|\partial_t U_1(x,l)|^2\nonumber\\
&&\3n\q\q\q\q\q  +|U_2(x,l)|^2+|\partial_t U_2(x,l)|^2\) dldx\\
&&\q\q +e^{2\zeta\phi_1}\frac{\lambda^2b
LM_0^2}{\pi} e^{\frac{\lambda^2}2
b^2}\(\max\{|\nabla\chi|^2, |\Delta\chi|^2\}
\int_{-\frac L2}^{\frac L2} \int_{O_3\setminus
O_2} \( |U_1(x,l)|^2 + |\nabla
U_1(x,l)|^2\nonumber\\
&&\3n\q\q\q\q\q + |U_2(x,l)|^2 + |\nabla
U_2(x,l)|^2\)dxd\nonumber\\
&&\3n\q\q +\frac{\lambda^2Lb}{2\pi}
\zeta^5\mu^7\phi_2^5e^{2\zeta\phi_2}
e^{\frac{\lambda^2}2 b^2}
\int_{\omega_0}\int_{-\frac L2}^{\frac L2}
|U_1(x,l)|^2dldx\]\nonumber\\
&&\3n\q\q  +Ce^{2\zeta\phi_4}
\int_{(-b,b_0)\cup(b_0,b)}\int_{\tilde\o}\big(|\partial_s
\wt U^F_1|^2+|\wt U^F_1|^2+|\partial_s \wt U^F_2|^2+|
\wt U^F_2|^2\big)dxds.\nonumber
\end{eqnarray}

Then
\begin{eqnarray}\label{925eq1}
&&\3n\3n\int_{-1}^{1} \int_{\tilde\omega}
\(| U^F_1|^2+| U^F_2|^2\)dxds\nonumber\\
&&\3n\3n\leq C\zeta^{-3}\mu^{-4}\phi_3^{-3}e^{-2\zeta\phi_3} \[e^{2\zeta\phi_2}\frac{4\lambda^2bL}{\pi} e^{\frac{\lambda^2}2 (b^2-(\frac L8)^2))} \max_{K}\{ |\Phi''(l)|^2, |\Phi'(l)|^2\}\nonumber\\
&&\3n\3n\q\q\q\q\q\q \times \int_{\tilde\omega} \int_{K} \(|U_1(x,l)|^2+|\partial_t U_1(x,l)|^2+|U_2(x,l)|^2+|\partial_t U_2(x,l)|^2\) dldx\nonumber\\
&&\3n\3n\q\q+e^{2\zeta\phi_1}\frac{\lambda^2b
LM_0^2}{\pi} e^{\frac{\lambda^2}2
b^2}\max\{|\nabla\chi|^2, |\Delta\chi|^2\}
\int_{-\frac L2}^{\frac L2} \int_{O_3\setminus
O_2} \( |U_1(x,l)|^2 + |\nabla
U_1(x,l)|^2\nonumber\\
&&\3n\3n\q\q+ |U_2(x,l)|^2 + |\nabla
U_2(x,l)|^2\)dxdl +\frac{\lambda^2Lb}{2\pi}
\zeta^5\mu^7\phi_2^5e^{2\zeta\phi_2}
e^{\frac{\lambda^2}2 b^2}
\int_{\omega_0}\int_{-\frac L2}^{\frac L2}
|U_1(x,l)|^2dldx\]\nonumber\\
&&\3n\3n\q +C\zeta^{-3}\mu^{-4}\phi_3^{-3}e^{-2\zeta(\phi_3-\phi_4)} \int_{(-b,b_0)\cup(b_0,b)}\int_{\tilde\o}\big(|\partial_s \wt U^F_1|^2+|\wt U^F_1|^2+|\partial_s \wt U^F_2|^2+|\wt U^F_2|^2\big)dxds\nonumber\\
&&\3n\3n\leq C\[\zeta^{-3}\mu^{-4}\phi_3^{-3}\(e^{2\zeta(\phi_2-\phi_3)}\frac{4\lambda^2bL}{\pi} e^{\frac{\lambda^2}2 (b^2-(\frac L8)^2)} \max_{K}\{ |\Phi''(l)|^2, |\Phi'(l)|^2\}+ e^{-2\zeta(\phi_3-\phi_4)}\nonumber\\
&&\3n\3n\q\q\q+e^{2\zeta(\phi_1-\phi_3)}\frac{\lambda^2b
LM_0^2}{\pi} e^{\frac{\lambda^2}2
b^2}\max\{|\nabla\chi|^2,
|\Delta\chi|^2\}\)C(M_0,M_1)^2\nonumber\\
&&\3n\3n\q\q\q+\zeta^2\mu^3\frac{\lambda^2Lb}{2\pi}
\frac{\phi_2^5}{\phi_3^3}e^{2\zeta(\phi_2-\phi_3)}
e^{\frac{\lambda^2}2 b^2}
\int_{\omega_0}\int_{-\frac L2}^{\frac L2}
|U_1(x,l)|^2dldx\].
\end{eqnarray}

Define
\begin{equation*}
\tau\triangleq
1-\exp\(-\mu\[\frac1M\(1-\frac{\beta_2}{\|\hat\psi\|_{L^\infty(\o)}}\)+1-b_0^2\]\),
\end{equation*}
then
\begin{equation}
\phi_3-\phi_4=\tau\phi_3>0.
\end{equation}
Let $\zeta=\frac{\lambda^2b^2}{4\tau\phi_3}$ and
$A>1$. By choosing $L=8Ab$, we have
\begin{eqnarray}\label{928eq1}
&&\int_{-1}^{1} \int_{\omega_1\setminus O_1}
\(| U^F_1|^2+| U^F_2|^2\)dxds\nonumber\\
&&\leq\int_{-1}^{1} \int_{\tilde\omega}
\(| U^F_1|^2+| U^F_2|^2\)dxds\nonumber\\
&&\3n\3n\leq C\[\zeta^{-3}\mu^{-4}\phi_3^{-3}\(e^{2\zeta(\phi_2-\phi_3)}\frac{32A\lambda^2b^2}{\pi} e^{-\frac{\lambda^2b^2}2(A^2-1)} \max_{K}\{ |\Phi''(l)|^2, |\Phi'(l)|^2\}+ e^{-2\zeta(\phi_3-\phi_4)}\nonumber\\
&&\3n\3n\q\q\q+e^{2\zeta(\phi_1-\phi_3)}\frac{8A\lambda^2b^2
M_0^2}{\pi} e^{\frac{\lambda^2}2
b^2}\max\{|\nabla\chi|^2,
|\Delta\chi|^2\}\)C(M_0,M_1)^2\nonumber\\
&&\3n\3n\q\q\q+\zeta^2\mu^3\frac{4A\lambda^2b^2}{\pi}
\frac{\phi_2^5}{\phi_3^3}e^{2\zeta(\phi_2-\phi_3)}
e^{\frac{\lambda^2}2 b^2}
\int_{\omega_0}\int_{-4Ab}^{4Ab}
|U_1(x,l)|^2dldx\] \nonumber\\
&&\3n\3n\leq C\[\(\frac{4\tau\phi_3}{\lambda^2b^2}\)^3\mu^{-4}\phi_3^{-3}\(e^{2\frac{\lambda^2b^2}{4\tau\phi_3}(\phi_2-\phi_3)}\frac{32A\lambda^2b^2}{\pi} e^{-\frac{\lambda^2b^2}2(A^2-1)} \max_{K}\{ |\Phi''(l)|^2, |\Phi'(l)|^2\}\nonumber\\
&&\3n\3n\q\q\q+ e^{-\frac{\l^2
b^2}{2}}+e^{2\frac{\lambda^2b^2}{4\tau\phi_3}(\phi_1-\phi_3)}\frac{8A\lambda^2b^2
M_0^2}{\pi} e^{\frac{\lambda^2}2
b^2}\max\{|\nabla\chi|^2,
|\Delta\chi|^2\}\)C(M_0,M_1)^2\nonumber\\
&&\3n\3n\q\q\q+\(\frac{\lambda^2b^2}{4\tau\phi_3}\)^2\mu^3\frac{4A\lambda^2b^2}{\pi}
\frac{\phi_2^5}{\phi_3^3}e^{2\frac{\lambda^2b^2}{4\tau\phi_3}(\phi_2-\phi_3)}
e^{\frac{\lambda^2}2 b^2}
\int_{\omega_0}\int_{-4Ab}^{4Ab}
|U_1(x,l)|^2dldx\] \nonumber\\
&&\3n\3n\leq C\[\(\frac{4\tau}{\lambda^2b^2}\)^3\mu^{-4}\(\frac{32A\lambda^2b^2}{\pi} e^{-(A^2-1-\frac{\phi_2-\phi_3}{\tau\phi_3})\frac{\lambda^2b^2}2} \max_{K}\{ |\Phi''(l)|^2, |\Phi'(l)|^2\}\nonumber\\
&&\3n\q\q\q+
e^{-\frac{\lambda^2b^2}{2}}+\frac{8A\lambda^2b^2
M_0^2}{\pi}e^{-\frac{\lambda^2b^2}{2}(\frac{\phi_3-\phi_1}{\tau\phi_3}-1)}\max\{|\nabla\chi|^2,
|\Delta\chi|^2\}\)C(M_0,M_1)^2\nonumber\\
&&\3n\3n\q\q\q+\frac{A\l^6b^6\mu^3}{4\pi\tau^2}\(\frac{\phi_2}{\phi_3}\)^5
e^{\frac{\lambda^2b^2}{2}(\frac{\phi_2-\phi_3}{\tau\phi_3}+1)}
\int_{\omega_0}\int_{-4Ab}^{4Ab}
|U_1(x,l)|^2dldx\]
\end{eqnarray}

By Parseval's identity, we get that
\begin{equation}\label{1004eq1}
\begin{array}{ll}
\ds\q\|\Phi U_j\|_{L^2((\omega_1\setminus O_1)\times(-\frac L2,\frac L2))}^2=\int_{-\frac L2}^{\frac L2}\int_{\omega_1\setminus O_1} |\Phi(t)U_j(x,t)|^2 dxdt\\
\ns\ds=\int_{\mathbb{R}}\int_{\omega_1\setminus O_1} |\Phi(t)U_j(x,t)|^2 dxdt=\frac{1}{2\pi}\int_{\mathbb{R}}\int_{\omega_1\setminus O_1} |\widehat{\Phi(l_0)U_j}(x,l_0)(t)|^2 dxdt\\
\ns\ds\leq\frac{1}{
\pi}\int_{\mathbb{R}}\int_{\omega_1\setminus O_1}
|(1-F_\lambda)\widehat{\Phi(l_0)U_j}(x,l_0)(t)|^2
dxdt+2\int_{\mathbb{R}}\int_{\omega_1\setminus O_1} |F_\lambda*\Phi(\cdot)U_j(x,\cdot)(l_0)|^2
dxdl_0.
\end{array}
\end{equation}
The first term in the right hand side of
\eqref{1004eq1} reads
\begin{eqnarray}\label{1004eq2}
&&\ds\frac{1}{\pi}\int_{\mathbb{R}}\int_{\omega_1\setminus O_1} |(1-F_\lambda)\widehat{\Phi(l_0)U_j}(x,l_0)(t)|^2 dxdt\nonumber\\
&&=\frac{1}{\pi}\int_{\mathbb{R}}\int_{\omega_1\setminus O_1} (1-e^{-(\frac t\lambda)^2})^2|\widehat{\Phi(l_0)U_j}(x,l_0)(t)|^2 dxdt\nonumber\\
&&\ds\leq\frac{2}{\pi\lambda^2}\int_{\mathbb{R}}\int_{\omega_1\setminus O_1} |t\widehat{\Phi(l_0)U_j}(x,l_0)(t)|^2 dxdt\nonumber\\
&&\ds\leq\frac{4}{\lambda^2}\int_{\mathbb{R}}\int_{\omega_1\setminus O_1} |\Phi'(l_0)u_j(x,l_0)+\Phi(l_0)\partial_{l_0} U_j(x,l_0)|^2 dxdl_0\\
&&\ds\leq\frac{8}{\lambda^2}\int_{\mathbb{R}}\int_{\omega_1\setminus O_1}\left( |\Phi'(l_0)U_j(x,l_0)|^2+\Phi(l_0)\partial_{l_0} U_j(x,l_0)|^2\right) dxdl_0\nonumber\\
&&\ds\leq\frac{8}{\lambda^2}\[\(\frac{1}{4L}\)^2\int_{K_0}\int_{\omega_1\setminus O_1}|U_j(x,l_0)|^2dxdl_0+\int_{-\frac L2}^{\frac L2}\int_{\omega_1\setminus O_1} |\partial_{l_0}U_j(x,l_0)|^2 dxdl_0\]\nonumber\\
&&\ds\leq\frac{8}{\lambda^2}\[\(\frac1{32Ab}\)^2\int_{K_0}\int_{\omega_1\setminus O_1}|U_j(x,l_0)|^2dxdl_0+\int_{-4Ab}^{4Ab}\int_{\omega_1\setminus O_1}
|\partial_{l_0}U_j(x,l_0)|^2 dxdl_0\].\nonumber
\end{eqnarray}

Let
\begin{equation}\label{512eq1}
U^F_{j,\lambda}(x,l_0)\triangleq
U^F_j(x,0)=\int_{\mathbb{R}}
F_\lambda(l_0-l)\Phi(l)U_j(x,l)dl=F_\lambda*\Phi(\cdot)U_j(x,\cdot)(l_0).
\end{equation}
By applying the Cauchy integral formula, for
$\rho\in(0,1)$ and by setting $z=\kappa+\rho
e^{i\phi}$, we have that
\begin{equation}\label{512eq2}
\begin{array}{ll}\ds
 U^F_{j,\lambda}(x,\kappa)\3n&\ds=\frac1{2\pi i}\int_{|z-\kappa|=\rho}\frac{U^F_{j,\lambda}(x,z)}{z-\kappa}dz \\
\ns&\ds=\frac1{2\pi i}\int_0^{2\pi}U^F_{j,\lambda}(x,\kappa+\rho e^{ir})dr \\
\ns&\ds=\frac1{2\pi i}\int_0^{1}\int_0^{2\pi}U^F_{j,\lambda}(x,\kappa+\rho e^{ir})dr d\rho\\
\ns&\ds=\frac1{2\pi
i}\int_{-1}^{1}\int_{-\sqrt{1-l_0^2}}^{\sqrt{1-l_0^2}}U^F_j(x,s)dsdl_0. 
\end{array}
\end{equation}
Thus,
\begin{equation}\label{1004eq3}
\begin{array}{ll}
\ds |U^F_{j,\lambda}(x,\kappa)|^2
\leq\frac{1}{\pi^2}\int_{-1}^{1}\int_{-1}^{1}\big|U^F_j(x,s)\big|^2dsdl_0.
\end{array}
\end{equation}
Integrating \eqref{1004eq3} with respect to $x$
over $\omega_1$, and with respect to $\kappa$ over
$[-\frac L2,\frac L2]$ we get that
\begin{equation}\label{1004eq4}
\begin{array}{ll}
\ds\int_{-4Ab}^{4Ab}\int_{\omega_1\setminus O_1} |U^F_{j,\lambda}(x,\kappa)|^2dxd\kappa\\
\ns\ds\leq\frac{1}{\pi^2} \int_{-4Ab}^{4Ab}\int_{-1}^{1}\left(\int_{-1}^{1}\int_{\omega_1\setminus O_1}\big|U^F_j(x,s)\big|^2dxds\right)dl_0d\kappa\\
\ns\ds\leq\frac{16Ab}{\pi^2}
\int_{-1}^{1}\int_{\omega_1\setminus O_1}\big|U^F_j(x,s)\big|^2dxds.
\end{array}
\end{equation}
Substituting \eqref{1004eq2} and \eqref{1004eq4}
into \eqref{1004eq1} yields
\begin{equation}\label{1001eq1}
\begin{array}{ll}
\ds\q\|\Phi U_j\|_{L^2((\omega_1\setminus O_1)\times(-4A,4A))}^2\\
\ns\ds\leq\frac{8}{\lambda^2}\[\frac1{(32Ab)^2}\int_{K_0}\int_{\omega_1\setminus O_1}|U_j(x,l_0)|^2dxdl_0+\int_{-4Ab}^{4Ab}\int_{\omega_1\setminus O_1} |\partial_{l_0}U_j(x,l_0)|^2 dxdl_0\]\\
\ns\ds\q\q+\frac{16Ab}{\pi^2}
\int_{-1}^{1}\int_{\omega_1\setminus O_1}\big|U^F_j(x,s)\big|^2dxds.
\end{array}
\end{equation}

Let $C_2=\max\{|\nabla\chi|^2, |\Delta\chi|^2\}$
and suppose that $A$ is  sufficiently enough
such that
$\frac12(A^2-1)\tau+1>\exp \mu\(\frac{\hat\psi(\hat
x)-\beta_2}{M\|\hat\psi\|_{L^\infty(\tilde\o)}}+1\)$.
Then
\begin{eqnarray}\label{1004eq5}
&&\|(\Phi U_1,\Phi U_2)\|^2_{L^2((\omega_1\setminus O_1)\times(-\frac L2,\frac L2))}\nonumber\\
&&\leq\frac{8}{\lambda^2}\sum_{j=1}^2\[\frac1{(32Ab)^2}\int_{K_0}\int_{\omega_1\setminus O_1}|U_j(x,l_0)|^2dxdl_0+\int_{-4Ab}^{4Ab}\int_{\omega_1\setminus O_1} |\partial_{l_0}U_j(x,l_0)|^2 dxdl_0\]\nonumber\\
&&\q\q+\sum_{j=1}^2\frac{16Ab}{\pi^2}
\int_{-1}^{1}\int_{\omega_1\setminus O_1}\big|U^F_j(x,s)\big|^2dxds\nonumber\\
&&\leq\frac{8}{\lambda^2}\sum_{j=1}^2\[\frac1{(32Ab)^2}\int_{K_0}\int_{\omega_1\setminus O_1}|U_j(x,l_0)|^2dxdl_0+\int_{-4Ab}^{4Ab}\int_{\omega_1\setminus O_1} |\partial_{l_0}U_j(x,l_0)|^2 dxdl_0\] \nonumber\\
&&\q\q+\frac{16Ab}{\pi^2}
C\Big\{\[\(\frac{4\tau}{\lambda^2b^2}\)^3\mu^{-4}\(\frac{32A\lambda^2b^2}{\pi} e^{-(A^2-1-\frac{\phi_2-\phi_3}{\tau\phi_3})\frac{\lambda^2b^2}2} \max_{K}\{ |\Phi''(l)|^2, |\Phi'(l)|^2\}\nonumber\\
&&\3n\q\q\q+
e^{-\frac{\lambda^2b^2}{2}}+\frac{8A\lambda^2b^2
M_0^2}{\pi}e^{-\frac{\lambda^2b^2}{2}(\frac{\phi_3-\phi_1}{\tau\phi_3}-1)}\max\{|\nabla\chi|^2,
|\Delta\chi|^2\}\)C(M_0,M_1)^2\nonumber\\
&&\3n\3n\q\q\q+\frac{A\l^6b^6\mu^3}{4\pi\tau^2}\(\frac{\phi_2}{\phi_3}\)^5
e^{\frac{\lambda^2b^2}{2}(\frac{\phi_2-\phi_3}{\tau\phi_3}+1)}
\int_{\omega_0}\int_{-4Ab}^{4Ab}
|U_1(x,l)|^2dldx\]\Big\}\nonumber\\
&&\leq\[\frac{8}{\lambda^2}
+C\frac{1024A\tau^3}{\mu^4\pi^2\lambda^6b^5}\(\frac{32A\lambda^2b^2}{\pi}e^{-\big(A^2-1-\frac{\phi_2-\phi_3}{\tau\phi_3}\big)\frac12\lambda^2b^2}+e^{-\frac12\lambda^2b^2}\nonumber\\
&&\qq\q+C_2\frac{8AM_0^2\lambda^2b^2}{\pi}e^{\frac{\lambda^2b^2}{2}(\frac{\phi_1-\phi_3}{\tau\phi_3}+1)}\)\]C(M_0,M_1)^2\nonumber\\
&&\qq+C\frac{A\mu^3\lambda^6b^6}{4\tau^2\pi}
\(\frac{\phi_2}{\phi_3}\)^5e^{\big(\frac{\phi_2-\phi_3}{\tau\phi_3}+1\big)\frac12\lambda^2
b^2} \|U_1\|^2_{L^2(\o_0\times(-\frac L2,\frac
L2))}.
\end{eqnarray}

Similarly, we have
\begin{eqnarray}\label{1004eq6}
&&\|(\Phi V_1,\Phi V_2)\|^2_{L^2((\omega_1\setminus O_1)\times(-\frac L2,\frac L2))}\nonumber\\
&&\leq\[\frac{8}{\lambda^2}
+C\frac{1024A\tau^3}{\mu^4\pi^2\lambda^6b^5}\(\frac{32A\lambda^2b^2}{\pi}e^{-\big(A^2-1-\frac{\phi_2-\phi_3}{\tau\phi_3}\big)\frac12\lambda^2b^2}+e^{-\frac12\lambda^2b^2}\nonumber\\
&&\qq\q+C_2\frac{8AM_0^2\lambda^2b^2}{\pi}e^{\frac{\lambda^2b^2}{2}(\frac{\phi_1-\phi_3}{\tau\phi_3}+1)}\)\]C(M_0,M_1)^2\nonumber\\
&&\qq+C\frac{A\mu^3\lambda^6b^6}{4\tau^2\pi}
\(\frac{\phi_2}{\phi_3}\)^5e^{\big(\frac{\phi_2-\phi_3}{\tau\phi_3}+1\big)\frac12\lambda^2
b^2} \|V_1\|^2_{L^2(\o_0\times(-\frac L2,\frac
L2))}\nonumber
\end{eqnarray}
and
\begin{eqnarray}\label{98eq5}
&&\3n\3n\|(c_{11}-\tilde{c}_{11},c_{22}-\tilde{c}_{22})\|_{L^2(\O)}
\leq C \sum_{j=1}^2\|\partial_t^j(w_1,w_2) \|_{L^2((\omega_1\setminus O_1)\times (0,T))} \nonumber\\
&&\3n\3n\leq C\left(\|(U_1,U_2) \|_{L^2((\omega_1\setminus O_1)\times (0,T))}+\|(V_1,V_2) \|_{L^2((\omega_1\setminus O_1)\times (0,T))}\right) \nonumber\\
&&\3n\3n\leq C\left(\|(\Phi U_1,\Phi U_2)\|_{L^2((\omega_1\setminus O_1)\times (0,T))}+\|(\Phi V_1,\Phi V_2)\|_{L^2((\omega_1\setminus O_1)\times  (0,T))}\right) \nonumber\\
&&\3n\3n\leq 2\[\frac{8}{\lambda^2}
+C\frac{1024A\tau^3}{\mu^4\pi^2\lambda^6b^5}\(\frac{32A\lambda^2b^2}{\pi}e^{-\big(A^2-1-\frac{\phi_2-\phi_3}{\tau\phi_3}\big)\frac12\lambda^2b^2}+e^{-\frac12\lambda^2b^2}\nonumber\\
&&\q+C_2\frac{8AM_0^2\lambda^2b^2}{\pi}e^{\frac{\lambda^2b^2}{2}(\frac{\phi_1-\phi_3}{\tau\phi_3}+1)}\)\]C(M_0,M_1)^2\nonumber\\
&&\q+C\frac{A\mu^3\lambda^6b^6}{4\tau^2\pi}
\(\frac{\phi_2}{\phi_3}\)^5e^{\big(\frac{\phi_2-\phi_3}{\tau\phi_3}+1\big)\frac12\lambda^2
b^2}
\(\|U_1\|^2_{L^2(\o_0\times(-\frac L2,\frac L2))}+\|V_1\|_{L^2(\o_0\times(-\frac L2,\frac L2))}\) \nonumber\\
&&\3n\3n\leq 2\[\frac{8}{\lambda^2}
+C\frac{1024A\tau^3}{\mu^4\pi^2\lambda^6b^5}\(\frac{32A\lambda^2b^2}{\pi}e^{-\big(A^2-1-\frac{\phi_2-\phi_3}{\tau\phi_3}\big)\frac12\lambda^2b^2}+e^{-\frac12\lambda^2b^2}\nonumber\\
&&\q+C_2\frac{8AM_0^2\lambda^2b^2}{\pi}e^{\frac{\lambda^2b^2}{2}(\frac{\phi_1-\phi_3}{\tau\phi_3}+1)}\)\]C(M_0,M_1)^2\nonumber\\
&&\q+C\frac{A\mu^3\lambda^6b^6}{4\tau^2\pi}
\(\frac{\phi_2}{\phi_3}\)^5e^{\big(\frac{\phi_2-\phi_3}{\tau\phi_3}+1\big)\frac12\lambda^2
b^2} \sum_{j=1}^2\|\partial_t^jw_1
\|_{L^2(Q_{\o_0})}.\nonumber
\end{eqnarray}
Let $\lambda\geq 0$ be such that
\begin{equation*}\label{519eq1}
\|(
c_{11}-\tilde{c}_{11},c_{22}-\tilde{c}_{22})\|_{L^2(\O)}
\leq
\frac{C_3}{\lambda^2}C(M_0,M_1)+e^{C_4{\lambda^2}}
\sum_{j=1}^2\|\partial_t^jw_1
\|_{L^2(Q_{\o_0})},
\end{equation*}
where $C_3$ and $C_4$ are two constants
independent of $\lambda$. By taking
$$
\lambda= \(\frac {\sum_{j=1}^2 |\ln
\|\partial_t^jw_1 \|_{L^2(Q_{\o_0})}|}
{C_4}\)^\frac12,
$$
if $\|\partial_t^jw_1 \|_{L^2(Q_{\o_0})}$ is
small enough, then
\begin{equation}\label{519eq2}
\begin{array}{ll}
\ds\|(a-\tilde a,
c_{11}-\tilde{c}_{11},c_{22}-\tilde{c}_{22})\|_{L^2(\O)}\\
\ns\ds
\leq \frac{C_3C_4}{ \sum_{j=2}^3|\ln \|\partial_t^jw_1 \|_{L^2(Q_{\o_0})}|}C(M_0,M_1)+\sum_{j=1}^2\|\partial_t^jw_1 \|_{L^2(Q_{\o_0})} \\
\ns \ds\leq C\( \sum_{j=1}^2|\ln
\|\partial_t^jw_1
\|_{L^2(Q_{\o_0})}|^{-1}C(M_0,M_1)+\sum_{j=1}^2\|\partial_t^jw_1 \|_{L^2(Q_{\o_0})}\) \\
\ns \ds\leq C\( \sum_{j=1}^2|\ln
\|\partial_t^jw_1
\|_{L^2(Q_{\o_0})}|^{-1}C(M_0,M_1)+\sum_{j=1}^2\|\partial_t^jw_1
\|_{L^2(Q_{\o_0})}\).
\end{array}
\end{equation}
Otherwise, there exists a constant $m>0$ such
that $\|\partial_t^jw_1 \|_{L^2(Q_{\o_0})}\geq
m$. Thus, by \eqref{57eq10} we have
\begin{equation}\label{519eq3}
\begin{array}{ll}\ds
\|(
c_{11}-\tilde{c}_{11},c_{22}-\tilde{c}_{22})\|_{L^2(\O)}
\leq C(M_0,M_1)\\
\ns\ds=\frac{C(M_0,M_1)}{m}m\leq
C\sum_{j=1}^2\|\partial_t^jw_1
\|_{L^2(Q_{\o_0})}\leq
C\sum_{j=1}^2\|\partial_t^jw_1
\|_{L^2(Q_{\o_0})}.
\end{array}
\end{equation}
%
%

\section*{acknowledgement}
The first author thanks the support of the NSFC
(No. 11501086,11971093), the Fundamental
Research Funds for the Central Universities (No.
ZYGX2019J094) and the Science  Strength
Promotion Programme of UESTC. The second author
is supported by Grant-in-Aid for Scientific
Research (S) 15H05740 of Japan Society for the
Promotion of Science, and prepared with the
support of the ``RUDN University Program 5-100".

\end{document}